\input ./style/arxiv-general.cfg
\documentclass[aos,MSNbibl,secthm,nameyear,seceqn,dvips]{arximspdf}
\makeatletter
   \@ifpackageloaded{graphicx}{}{\usepackage{graphicx}}
\makeatother
\usepackage{mathrsfs}

\doi{10.1214/15-AOS1374}
\volume{44}
\issue{2}
\pubyear{2016}
\firstpage{515}
\lastpage{539}
\docsubty{FLA}

\makeatletter
\newproclaim{exa}{Example}
\def\half{\frac{1}{2}}
\def\halff{\tfrac{1}{2}}
\newtheorem{tm}{Theorem}
\newtheorem{cy}{Corollary}
\newtheorem{pn}{Proposition}
\newcommand{\bbeta}{\bolds{\beta}}
\newcommand{\bX}{\mathbf{X}}
\makeatother

\begin{document}
\begin{frontmatter}

\title{Local independence feature screening for nonparametric and
semiparametric models by~marginal empirical likelihood}
\runtitle{Nonparametric EL screening}

\begin{aug}
\author[A]{\fnms{Jinyuan}~\snm{Chang}\thanksref{M1,M2,T1}\ead[label=e1]{jinyuan.chang@unimelb.edu.au}},
\author[B]{\fnms{Cheng Yong}~\snm{Tang}\thanksref{M3,T2}\ead[label=e2]{yongtang@temple.edu}}
\and
\author[C]{\fnms{Yichao}~\snm{Wu}\thanksref{M4,T3}\corref{}\ead[label=e3]{wu@stat.ncsu.edu}}
\runauthor{J. Chang, C. Y. Tang and Y. Wu}
\affiliation{Southwestern University of Finance and Economics\thanksmark{M1}, University of Melbourne\thanksmark{M2},
Temple University\thanksmark{M3} and North
Carolina State University\thanksmark{M4}}
\address[A]{J. Chang\\
School of Statistics\\
Southwestern University of Finance\\
\quad and Economics\\
Chengdu, Sichuan 611130\\
China \\
and\\
School of Mathematics and Statistics\\
University of Melbourne\\
Parkville, Victoria 3010\\
Australia \\
\printead{e1}}
\address[B]{C. Y. Tang\\
Department of Statistics\\
Temple University\\
1810 North 13th Street\\
Philadelphia, Pennsylvania 19122-6083\\
USA\\
\printead{e2}}
\address[C]{Y. Wu\\
Department of Statistics\\
North Carolina State University\\
2311 Stinson Drive\\
Raleigh, North Carolina 27695-8203\\
USA\\
\printead{e3}}
\end{aug}
\thankstext{T1}{Supported by the
Fundamental Research Funds for the Central Universities, NSFC (Grant No. 11501462)
and a grant from the Australian Research Council. The bulk of this work was carried
out during Jinyuan Chang's Ph.D. study at Peking University.}
\thankstext{T2}{Supported in part by NSF Grants IIS-1546087 and SES-1533956, and the Fox School of Business and Management, Temple University.}
\thankstext{T3}{Supported in part by NSF Grant DMS-10-55210 and NIH/NCI
Grants R01-CA-149569 and P01-CA-142538.}

%
\received{\smonth{2} \syear{2015}}
%
\revised{\smonth{8} \syear{2015}}

%
\begin{abstract}
We consider an independence feature screening technique for identifying
explanatory variables that locally contribute to the
response variable in high-dimensional regression analysis. Without
requiring a specific parametric form of the underlying data model, our
approach accommodates
a wide spectrum of nonparametric and semiparametric model families.
To detect the local contributions of explanatory variables, our
approach constructs empirical likelihood locally in conjunction with
marginal nonparametric regressions.
Since our approach actually requires no estimation,
it is advantageous in scenarios such as the single-index models where
even specification and identification of a marginal model is an issue.
By automatically incorporating the level of variation of the
nonparametric regression and directly assessing the strength of data
evidence supporting
local contribution from each explanatory variable, our approach
provides a unique perspective for solving feature screening problems.
Theoretical analysis shows that our approach can handle data
dimensionality growing exponentially with the sample size.
With extensive theoretical illustrations and numerical examples, we
show that the local independence screening approach performs promisingly.
\end{abstract}

%
\begin{keyword}[class=AMS]
\kwd[Primary ]{62G99}
\kwd[; secondary ]{62H99}
\end{keyword}
\begin{keyword}
\kwd{Empirical likelihood}
\kwd{high-dimensional data analysis}
\kwd{nonparametric and semiparametric models}
\kwd{sure independence screening}
\end{keyword}
\end{frontmatter}

\section{Introduction}\label{sec1}

High-dimensional data are becoming increasingly available, and they
have triggered surging investigation and development of new theory and
methods; see \citet{Hastieetal2009book}, \citet{FanLv2009Sinica}
and \citet{Buhlmannvan2011} for overview and discussions.
Independence feature screening is a class of rapidly developing
approaches that is particularly useful in preliminary analysis for
preprocessing data to reduce the scale of high-dimensional statistical
problems; see, among others, \citet{FanLv2008JRSSB} and \citet{FanSong2010AOS} for the independence screening methods for linear and
generalized linear models, \citet{MaiZou2013Bioka} for variable screening in
classification problems, \citet{Zhuetal2011JASA}, \citet{Lietal2012AOS} and \citet{Lietal2012JASA} for feature screening methods using more general
types of correlations.

Broadly speaking, independence feature screening methods rely on
ranking estimations measuring the marginal contributions of explanatory
variables.
For example, \citet{FanLv2008JRSSB} and \citet{FanSong2010AOS} consider ranking
magnitudes of marginal estimators under some parametric models. In
linear and
generalized linear models, marginal estimator-based ranking can be
viewed as equivalent to marginal correlation based ranking
[\citet{FanLv2008JRSSB}, \citeauthor{ChangTangWu2013} (\citeyear{ChangTangWu2013})].
Various generalized versions of the correlation and conditional
correlation are also considered as ranking criteria
[\citet{Zhuetal2011JASA}, \citet{Lietal2012AOS}, \citet{Liuetal2014,Fanet2014}]. Recently, \citet{FangFengSongJASA2011} consider a
ranking measure based on aggregating local contributions from an
explanatory variable in a framework
of nonparametric additive models using marginal penalized splines approach.
In classification problems, \citet{MaiZou2013Bioka} propose to use the
so-called Kolmogrov filter to construct a ranking criterion based on aggregating
sample distributional discrepancies between the two groups of interest
at all observed values of a predictor.

We consider in this paper an independence feature screening method for
a general class of regression problems covering the nonparametric
additive models,
semiparametric single-index models and multiple-index models, and
varying coefficient models as special cases. There are two building
blocks for constructing the
screening criterion in our approach. The first one is the nonparametric
regression applied marginally on one explanatory variable at a time.
For overview of nonparametric regression methods, see \citet{fan1996} and
\citet{hardle1990}. The second building block is empirical likelihood
[\citet{Owen1988}]  constructed locally for the marginal nonparametric regression.
Instead of acquiring some marginal estimators, our approach is capable
of objectively and conveniently assessing the strength of data evidence for
testing the local contributions of a given explanatory variable.
Moreover, as  has been noted in the literature, an independence
feature screening procedure may miss explanatory variables that are
marginally unrelated
but jointly related to the response [\citet{FanLv2008JRSSB}]. To address
this issue, many iterative versions of feature screening methods have
been proposed. We borrow the idea
of \citet{Zhuetal2011JASA} and propose an iterative version of our local
independence feature screening procedure.

Our study carries innovative contributions from a few aspects.
First and foremost, the perspective of our approach is unique compared
with other existing ones. Our approach directly targets at quantifying
the strength of data evidence against the null hypothesis that
explanatory variables are not locally contributing to the response
variable. Hence, it actually requires no estimation. Moreover, as shown
in our theoretical analysis, the fundamental statistic in our approach
is self-Studentized, automatically incorporating variance of the
marginal statistical approach.
All existing approaches for nonparametric and semiparametric models
require estimating marginal contributions and incorporate no effect
from the level of variations of the marginal estimators. As a
consequence, ranking the nonstandardized magnitudes of the marginal
estimators may not best reflect the marginal contributions from
predictors. Additionally, there may be difficulties when identifying
the marginal effect becomes an issue, for example, in single- and
multiple-index models.
We show in our numerical examples that our approach outperforms others
especially when the signal of the marginal contribution is weak, and
when the variation of the response variable is more complex, all thanks
to the unique perspective of the proposed marginal empirical likelihood
approach.
Second, existing approaches are typically investigated within specific
families of models while our approach targets at detecting generic
local contributions to the response variable from an explanatory variable.
Thus, our method is suitable for capturing more general nonlinear
effects in explanatory variables for solving a broad range of
high-dimensional problems.
Our theoretical analysis establishes the validity for feature
screening in a general and broad setup, allowing data dimensionality to
grow exponentially with the sample size, and our numerical and real
data examples demonstrate that our method performs very promisingly.

Our investigation also contributes in solving challenging empirical
likelihood problems for high-dimensional nonparametric and
semiparametric statistical problems.
In existing literature, much effort has been devoted into extending
the empirical likelihood of Owen (\citeyear{Owen1988}, \citeyear{Owen1991}) for parametric models to
nonparametric and semiparametric models; see, among others, \citet{Chen1996}, \citet{ChenQin2000} and the review in \citet{ChenVan2009test}.
For high-dimensional data, it remains open for solving empirical
likelihood problems in nonparametric and semiparametric scenarios where
merits such as robustness and other nonparametric features are highly
desirable [\citet{Hjortetal2008AS,ChenPengQin2008,TangLeng2010Bioka,LengTang2010,ChangChenChen2012}].
Recently, \citeauthor{ChangTangWu2013} (\citeyear{ChangTangWu2013}) investigate marginal empirical
likelihood
for general high-dimensional parametric models
specified by the estimating equations.
Nevertheless, studying high-dimensional empirical
likelihood beyond parametric models remains open because formulating
and characterizing the empirical likelihood locally itself is known to
be an important and difficult problem [\citet{ChenVan2009test}].
Our study on the local feature screening procedure solves the problem
of constructing and characterizing empirical likelihood locally, which
ideally fits a broad class of nonparametric and semiparametric models.
Additionally, for summarizing the contribution from one specific
predictor, we propose and justify an approach for aggregating the data
evidence for local contributions, which in turn delivers the validity
of our feature screening procedure. Remarkably, our approach can handle
exponential data dimensionality even in the nonparametric and
semiparametric settings where the convergence rate of nonparametric
kernel regression is known to be slower [\citet{fan1996}, \citet{hardle1990}].

The rest of the paper is organized as follows. Section~\ref{s2} introduces the
methodology of independence feature screening for nonparametric models
and presents the corresponding theoretical properties. In
Section~\ref{sec3}, we apply this unified screening approach to deal with
problems in nonparametric additive models, single-index models and
multiple-index models, and
varying coefficient models.
As a methodological extension, we outline the iterative version of our
unified screening
approach in Section~\ref{sec4}. Our
simulation studies in Section~\ref{sec5} demonstrate the effectiveness of this
method. We conclude with a discussion in Section~\ref{s6}, and relegate a real
data analysis and the proofs to the supplementary file of this paper
[\citet{ChangTangWu2015}].

\section{Main results} \label{s2}
\subsection{Methods}\label{sec21}

Suppose that we have a random sample $\{(\bX_i,Y_i)\}_{i=1}^n$ from the
data model\vspace*{-3pt}
%
\begin{equation}
\label{eqm}
Y=m(\bX)+\varepsilon,
\end{equation}
where $\bX=(X_1,\ldots,X_p)^{\mathrm{T}}$ and $\varepsilon$ is the random error
with $\mathbb{E}(\varepsilon|\bX)=0$. In our study, no specification of
$m(\bX)$ is required.
The
data dimensionality $p$ of the explanatory variable vector $\bX$ can
grow exponentially with the sample size $n$, but the true model is very
sparse in the sense that there is only a small fraction of the
explanatory variables contributing to the response variable $Y$. Let
$
\mathcal{M}_*= \{1\leq j\leq p: \mathbb{E}(Y|\bX)$  varies
with the value of $X_j\}$,
and we call the variables indexed by $\mathcal{M}_*$ contributing
explanatory variables. Without loss of generality, we assume that
$\mathbb{E}(Y)=0$ implying that $\mathbb{E}\{m(\bX)\}=0$.

Since $\bX$ is high-dimensional and without
any prior information on which of them are contributing in explaining
$Y$, a natural idea is to investigate the marginal contribution from
each explanatory variable in explaining $Y$ to justify whether it is relevant.
For
such a purpose, we consider
marginal nonparametric regression problems:\vspace*{-3pt}
%
\begin{equation}\label{eqmj}
\min_{f_j\in\mathscr{L}_2} \mathbb{E}\bigl[\bigl\{Y-f_j(X_j)
\bigr\}^2\bigr] \qquad (j=1,\ldots, p),
\end{equation}
where $\mathscr{L}_2$ denotes the class of square integrable functions.
Note that $ \mathbb{E}(Y|X_j)$ is the minimizer of (\ref{eqmj}).
Naturally, we use ${f}_j(x)=\mathbb{E}(Y|X_j=x)$ to evaluate the
marginal contribution of $X_j$ locally at $X_j=x$.
If an explanatory variable $X_j$ is not contributing to $Y$
marginally, then ${f}_j(x)=0$
for all $x\in\mathcal{X}$. Here, $\mathcal{X}$ is the support of $X_j$.
This motivates us to investigate a feature screening procedure by
assessing whether $\mathbb{E}(Y|X_j)\equiv0$ or not for each $j=1,\ldots
,p$. However, how to develop a nice way for summarizing the impact due
to $X_j$ is not straightforward due to the fact that $ f_j(x)$ is a
function of $x$ so that one needs to assess $f_j(x)$ for all values
over the support of $X_j$.

We consider the Nadaraya--Watson (NW) estimator for\vspace*{-2pt} $f_j(x)$:
%
\begin{equation}\label{eqkernel}
\hat{f}_j(x)=\frac{n^{-1}\sum_{i=1}^n
\mathcal{K}_h(X_{ij}-x)Y_i}{n^{-1}\sum_{i=1}^n \mathcal
{K}_h(X_{ij}-x)},
\end{equation}
where $\mathcal{K}_h(u)=h^{-1}\mathcal{K}(\frac{u}{h})$ for some kernel
function $\mathcal{K}$ and $h$ is the bandwidth. For capturing the
marginal variable effect, the choice of the NW estimator does not
compromise the general applicability of the marginal empirical
likelihood with other nonparametric approaches,
for example, the local linear estimator [\citet{fan1996}], etc.
Intuitively, $\hat{f}_j(x)$ should be small for all $x\in\mathcal{X}$ if
$X_j$ does not marginally contribute to $Y$.

Empirical likelihood [Owen (\citeyear{Owen1988}, \citeyear{Owen2001})] is an influential nonparametric
likelihood approach. Without requiring to assume full parametric
distributions, empirical likelihood shares some desirable merits of the
conventional likelihood such as $\chi^2$-distributed likelihood ratios
and Bartlett correctability; see
\citet{ChenVan2009test} for a review.
For assessing $f_j(x)=0$ at a given $x$ without distributional
assumptions, we construct the following empirical\vspace*{-2pt} likelihood:
%
\begin{eqnarray}
&&\mathrm{EL}_j(x,0)
\nonumber
\\[-9pt]
\label{local}
\\[-9pt]
\nonumber
&&\qquad =\sup \Biggl\{\prod_{i=1}^nw_i:w_i
\geq0, \sum_{i=1}^nw_i=1,\sum
_{i=1}^nw_i
\mathcal{K}_h(X_{ij}-x)Y_i=0 \Biggr\}.
\end{eqnarray}
By applying the Lagrange multiplier method for solving (\ref{local}),
we obtain the empirical likelihood\vspace*{-2pt} ratio:
%
\begin{eqnarray}
\ell_j(x,0) &=& -2\log\bigl\{\mathrm{EL}_j(x,0)\bigr
\}-2n\log n
\nonumber
\\[-9pt]
\label{eqmelr}
\\[-9pt]
\nonumber
&=& 2\sum_{i=1}^n\log\bigl\{1+\lambda
\mathcal{K}_h(X_{ij}-x)Y_i\bigr\},
\end{eqnarray}
where $\lambda$ is the univariate Lagrange multiplier solving $ 0=\sum_{i=1}^n\frac{\mathcal{K}_h(X_{ij}-x)Y_i}{1+\lambda \mathcal
{K}_h(X_{ij}-x)Y_i}$.
Since the denominator in (\ref{eqkernel}) converges to the density of
$X_j$ evaluated at $x$, a~large value of $\ell_j(x,0)$ is taken as
evidence against $f_j(x)=0$ provided that the density of $X_j$ is
bounded away from 0 at $x$.
Hence, $\ell_j(x,0)$ is indeed a statistic for testing whether or not
the numerator in (\ref{eqkernel}) has zero mean locally at $x$.
If $0$ is not in the convex hull of $\{\mathcal{K}_h(X_{ij}-x)Y_i\}
_{i=1}^n$, we define $\ell_j(x,0)=\infty$ as a strong evidence of the
local contribution from $X_j$.

For assessing $\mathbb{E}(Y|X_j)\equiv 0$, we propose to\vspace*{-2pt} use
%
\begin{eqnarray}\label{eqel1}
\ell_j(0)&=& \sup_{x\in\mathcal{X}_n}\ell_j(x,0)
\end{eqnarray}
for each $j=1,\ldots,p$, where $\mathcal{X}_n$
is a partition of the support $\mathcal{X}$ into several intervals.
For feature screening purpose, we propose selecting the set of
explanatory variables by\vspace*{-2pt}
%
\begin{equation}
\label{eqselect}
\widehat{\mathcal{M}}_{\gamma_n}= \bigl\{1\leq j\leq p:
\ell_j(0)\geq \gamma_n \bigr\}.
\end{equation}
We call our method a local independence feature screening approach by
observing that in (\ref{local}) it is the local correlation evaluated
at $x$ between $X_j$ via $\mathcal{K}_h(X_j-x)$ and $Y$ is assessed.
Such a notion is seen as extending the initial proposal of independence
feature screening using correlations as in \citet{FanLv2008JRSSB}. As the
empirical likelihood ratio $\ell_j(x,0)$ is self-Studentized, the
statistic $\ell_j(0)$ is robust to the heterogeneous cases. Compared
with $L_2$-type screening statistics,
the $L_\infty$-type screening statistic $\ell_j(0)$ is seen to be more suitable
when $\mathbb{E}(Y|X_j=x)$ is away from zero
locally at some $x$ instead of globally. We will give the specification
of $\gamma_n$ later in Theorem~\ref{tm1} under which the proposed approach has
the sure screening property, that is, capable of identifying the set of
contributing explanatory variables. However, choosing
$\gamma_n$ is generally hard in practice. Thus, following the
convention of existing screening methods, we suggest running
procedure to preprocess data by selecting a prespecified number of
variables.

Practically implementing the proposed method is convenient. In
practice, a natural choice is to evaluate the statistic\vspace*{1pt} (\ref{eqel1})
by $\ell_j(0)=\max_{1\leq i \leq n} \ell_j(X_{ij},0)$, where $\{X_{ij}\}
_{i=1}^n$ are the $n$ observations of the $j$th explanatory variable.
Since evaluating $\ell_j(x,0)$ only involves univariate optimization
when solving (\ref{eqmelr}) using the Lagrange multiplier method, the
screening statistics can be carried out easily by the existing
algorithms. As for the bandwidth $h$ in the marginal NW estimator (\ref
{eqkernel}), we note that conventional bandwidth selection methods such
as cross-validation and the reference rule [\citet{fan1996}] can
be applied.
Our theory in the next section demonstrates the validity of the
variable screening procedure for a range of bandwidth applied in (\ref
{eqkernel}) including ones selected by methods like cross-validation
and the reference rule. Our numerical examples
also show that the approach implemented with bandwidth selected by the
cross-validation performs satisfactorily.

\subsection{Theoretical properties}\label{the}

Throughout this paper, we use $\|\cdot\|_2$ and $\|\cdot\|_\infty$ to
denote the $L_2$-norm and sup-norm, respectively, and $C^r(\mathcal
{I})$ denotes the class of all continuous functions defined over
$\mathcal{I}$ that are $r$ times differentiable.
We assume the following conditions.
\begin{longlist}[(A.6)]
\item[(A.1)]  The
marginal projections $\{f_j\}_{j=1}^p$ belong to $C^r(\mathcal{X})$. If $r=0$,
$f_j$'s satisfy the Lipschitz condition with order $\alpha\in(0,1]$,
that is $ |f_j(s)-f_j(t)|\leq K_1|s-t|^\alpha$ for any $s,t\in
\mathcal{X}$, where $K_1$ is a positive constant uniformly for any
$j=1,\ldots,p$. In addition, there exists a constant $K_2$ such that
$|f_j^{(r)}(x)|\leq K_2$ for any $x\in\mathcal{X}$ and $j=1,\ldots,p$.

\item[(A.2)] The marginal density function $g_j$ of $X_j$ satisfies
$0<K_3\leq g_j(x)\leq K_4<\infty$ on $\mathcal{X}$ for $j=1,\ldots,p$.
In addition, we assume that each $g_j$ belongs to $C^r(\mathcal{X})$
for the $r$ given in (A.1) and $|g_j^{(r)}(x)|\leq K_5$ for any $x\in
\mathcal{X}$ and $j=1,\ldots,p$.

\item[(A.3)]  There exist nonnegative constants $c_1>0$ and $\kappa\in[0,\frac
{\max\{r,\alpha\}}{2\max\{r,\alpha\}+2})$ such that $ \min_{j\in\mathcal
{M}_*}\|f_j\|_\infty\geq c_1n^{-\kappa}$, where $r$ and $\alpha$ are
specified in (A.1).

\item[(A.4)]  Let $\|\mathcal{X}_n\|$ be the largest length of the intervals
in the partition $\mathcal{X}_n$, and there exists some positive
constant $\xi$ such that $\|\mathcal{X}_n\|=n^{-\xi}$.

\item[(A.5)] There exist positive constants $K_6$, $K_7$ and $\gamma$ such
that $
\mathbb{P}(|Y|\geq u)\leq K_6\exp(-K_7u^\gamma)
$ for any $u>0$.

\item[(A.6)]  For $r$ specified in (A.1), if $r\geq1$, the kernel function
$\mathcal{K}(\cdot)$ is of order $r$, that is, $ \int \mathcal{K}(u)\,du=1$,
$\int u^k\mathcal{K}(u)\,du=0$ for $k=1,\ldots,r-1$ and $\int u^r\mathcal
{K}(u)\,du>0$. If $r=0$, the kernel function satisfies $\mathcal{K}(u)\geq
0$ for any $u$ and $\int \mathcal{K}(u)\,du=1$.
\end{longlist}

Here, (A.1) is a general condition describing the continuity of each
$f_j(x)=\mathbb{E}(Y|X_j=x)$. If the first derivation of each $f_j$
exists, then $r\geq1$ and $\alpha=1$. Assumption (A.2) is standard for
kernel regression implying that the density of $X_j$ does not vanish on
its support; see, for example, \citet{hardle1990}, and it implies
bounded support of the explanatory variables.
For ease in presentation,
we take the same support $\mathcal{X}=[a,b]$ for all $X_j$ which can be
easily satisfied because some location-scale transformation can always
be applied in practice if otherwise. The condition in (A.3) is for
identifying $\mathcal{M}_*$, which requires that the minimal signal
strength measured by $\|f_j\|_\infty$ cannot be too weak. The
restriction of the minimal signal strength depends on the continuity of
$f_j$ via $r$.
The smoother $f_j$'s are, the weaker the condition on the signal
strength is required, and the minimal signal strength cannot vanish at
a rate faster than $n^{-{1}/{2}}$. Assumption (A.4) regularizes the
partition of the support $\mathcal{X}$ to be of size at least $O(n^\xi)$.
Assumption (A.5) on the tail distribution of the response variable is a
conventional technical requirement for Cram\'{e}r-type large
deviations. For example, $\gamma=2$ if the response variable $Y$ is a
normal or sub-Gaussian distribution, and $\gamma=\infty$ if $Y$ has a
compact support. Assumption (A.6) specifies the requirement for the
kernel function so that the bias due to kernel smoothing is not
dominating; see \citet{Muller1987} for more detail about higher order
kernel functions.

For the parameters $r, \alpha$ and $\gamma$ specified in above
assumptions, let
%
\begin{equation}
\label{eqparameters}
\qquad\varrho_1=\max\{r,\alpha\}, \qquad\varrho_2=
\max\{\gamma,2\}+2 \quad\mbox{and}\quad \delta=\max\biggl\{\frac{2}{\gamma}-1,0\biggr\}.
\end{equation}
The parameter $\varrho_1$ characterizes the continuity of the marginal
projections and densities. Parameters
$\varrho_2$ and $\delta$ are related to the tail probabilistic
behavior of response variable $Y$. Meanwhile, we assume that the
bandwidth $h$ used in (\ref{local}) satisfies
$
h\asymp n^{-w}
$
for some positive $w$ whose specification is discussed later.

\begin{pn}\label{pn1}
Under assumptions \textup{(A.1)}--\textup{(A.6)}, pick $w\in[\frac{\kappa
}{\varrho_1},1)$ and $\xi>\kappa+2w$, then there exists a uniform
constant $C_1$ such that for any $j\in\mathcal{M}_*$ and $L\rightarrow
\infty$,
\begin{eqnarray*}
\mathbb{P} \biggl\{\ell_j(0)<\frac{c_1^2K_3^2n^{1-2\kappa-2w}}{2L^2} \biggr\} &\leq &  \exp
\bigl(-C_1L^\gamma \bigr)
\\
&&{}+\exp \bigl(-C_1n^{\min \{1-2\kappa-w,({1-\kappa-w})/({1+\delta})\}} \bigr).
\end{eqnarray*}
\end{pn}

Proposition~\ref{pn1} gives a uniform result for all explanatory
variables contributing in the true model. The maximum distance between
the adjacent two points that are used to construct our procedure should
be $o(n^{-\kappa-2w})$. Specifically, with large probability and
uniformly for all $j\in \mathcal{M}_*$, the diverging rate of $\ell
_j(0)$ is not slower than $n^{1-2\kappa-2w}L^{-2}$.
If $j\notin \mathcal{M}_*$, that is, the explanatory variable $X_j$
does not have the marginal contribution to $Y$ (i.e., $f_j=0$),
following the argument in \citet{Owen1988} and \citeauthor{ChangTangWu2013} (\citeyear{ChangTangWu2013}),
it can be shown that the corresponding $\ell_j(0)$ is $O_p(1)$. Hence,
$n^{{1}/{2}-\kappa-w}L^{-1}$ is required to diverge as $n\rightarrow
\infty$ for sure independence screening. Furthermore, we note that the
requirement for the bandwidth used in Proposition~\ref{pn1} is mild, which can
be naturally satisfied by the conventional optimal bandwidth
$h=O(n^{-{1}/{5}})$ selected by the cross-validation method.

Let $L=n^{{1}/{2}-\kappa-w-\tau}$ for some $\tau\in(0,\frac
{1}{2}-\kappa-w)$.
A more clear uniform result
related to the probabilistic behavior of the statistics $\ell_j(0)$
for $j\in\mathcal{M}_*$ is described in the following corollary.

\begin{cy}\label{cy1}
Under assumptions \textup{(A.1)}--\textup{(A.6)}, pick $w\in[\frac{\kappa
}{\varrho_1},\frac{1}{2}-\kappa)$, $\tau\in(0,\frac{1}{2}-\kappa-w)$ and
$\xi>\kappa+2w$, then
\begin{eqnarray*}
\max_{j\in\mathcal{M}_*}\mathbb{P} \bigl\{\ell_j(0)<\halff
c_1^2K_3^2n^{2\tau} \bigr
\}& \leq & \exp \bigl\{-C_1n^{({1}/{2}-\kappa
-w-\tau)\gamma} \bigr\}
\\
&&{}+\exp \bigl(-C_1n^{\min \{1-2\kappa-w,({1-\kappa-w})/({1+\delta})
\}} \bigr),
\end{eqnarray*}
where $C_1$ is given in Proposition~\ref{pn1}.
\end{cy}

Choosing threshold level $\gamma_n=\frac{1}{2}c_1^2K_3^2n^{2\tau}$ in
(\ref{eqselect}) and noting that
\begin{eqnarray*}
\mathbb{P} (\mathcal{M}_*\varsubsetneq\widehat{\mathcal{M}}_{\gamma
_n} ) &=&
\mathbb{P} \bigl\{\mbox{There exists } j\in\mathcal {M}_* \mbox{ such that }
\ell_j(0)<\halff c_1^2K_3^2n^{2\tau}
\bigr\}
\\
&\leq &  s\max_{j\in\mathcal{M}_*}\mathbb{P} \bigl\{
\ell_j(0)< \halff c_1^2K_3^2n^{2\tau}
\bigr\},
\end{eqnarray*}
we establish the sure screening property of our approach in the
following theorem based on Corollary~\ref{cy1}.

\begin{tm}
\label{tm1}
Under assumptions \textup{(A.1)}--\textup{(A.6)}, pick $w\in[\frac{\kappa
}{\varrho_1},\frac{1}{2}-\kappa)$, $\gamma_n=\frac
{1}{2}c_1^2K_3^2n^{2\tau}$ for some $\tau\in(0,\frac{1}{2}-\kappa-w)$, and
$\xi>\kappa+2w$, then
\begin{eqnarray*}
\mathbb{P} (\mathcal{M}_*\subset\widehat{\mathcal{M}}_{\gamma
_n} ) &\geq & 1
-s\exp \bigl\{-C_1n^{({1}/{2}-\kappa-w-\tau)\gamma
} \bigr\}
\\
&&{}-s\exp \bigl(-C_1n^{\min \{1-2\kappa-w,({1-\kappa-w})/({1+\delta})
\}} \bigr),
\end{eqnarray*}
where $C_1$ is given in Proposition~\ref{pn1}.
\end{tm}

Theorem~\ref{tm1} implies
that our local independence feature screening method can handle
nonpolynomial dimensionality:
$
\log p
=o(n^{\epsilon})$
for $\epsilon=\min \{1-2\kappa-w,(\frac{1}{2}-\kappa-w-\tau)\gamma
 \}$.
By noting that $w\geq \frac{\kappa}{\varrho_1}$, the highest
dimensionality is achieved with the optimal $\epsilon=\min\{1-2\kappa
-\frac{\kappa}{\varrho_1},(\frac{1}{2}-\kappa-\frac{\kappa}{\varrho
_1})\gamma\}$ when $\tau$ is close enough to zero. It actually depends
on $\kappa$, $\varrho_1$ and $\gamma$, that is, the signal strength,
smoothness of $f_j$'s and the tail probabilistic behavior of $Y$. If
$Y$ follows a normal or sub-Gaussian distribution such that $\gamma=2$,
the corresponding highest dimensionality satisfies $\log
p=o(n^{1-2\kappa-{2\kappa}/{\varrho_1}})$. Furthermore, if the
projections $f_j$'s have derivatives of all orders such that $\varrho
_1=r=\infty$, then the highest dimensionality satisfies $\log
p=o(n^{1-2\kappa})$.

In what follows, we consider the size of the selected set $\widehat
{\mathcal{M}}_{\gamma_n}$ under an ideal case that
%
\begin{equation}\label{eqfa}
\max_{j\notin\mathcal{M}_*}\|f_j\|_\infty=o
\bigl(n^{-\kappa}\bigr).
\end{equation}
The key is to investigate the probabilistic behavior of $\mathbb{P}\{
\ell_j(0)\geq \half c_1^2K_3^2n^{2\tau}\}$ for each $j\notin\mathcal
{M}_*$ which is given in the next proposition.

\begin{pn}\label{pn2}
Under assumptions \textup{(A.1)}--\textup{(A.2)} and \textup{(A.4)}--\textup{(A.6)}, suppose $\max_{j\notin\mathcal{M}_*}\|f_j\|_\infty=O(n^{-\eta})$ for some $\eta>\frac
{(2\varrho_1+1)\kappa}{2\varrho_1}$. Pick $w\in[\frac{\kappa}{\varrho
_1},\break \min\{\frac{1}{2}-\kappa,2(\eta-\kappa)\})$, $\tau\in(\max\{\frac
{1}{2}-\eta-\frac{w}{2},0\},\frac{1}{2}-\kappa-w)$ and $\xi>\kappa+2w$.
If $\inf_{u\in[a,b]}\mathbb{E}(Y^2|X_j=u)\geq\rho$ for some positive
$\rho$ holds for any $j\notin\mathcal{M}_*$, then there exists a
uniform positive constant $C_2$ such that for any $j\notin\mathcal{M}_*$,
\[
\mathbb{P} \bigl\{\ell_j(0)\geq \halff c_1^2K_3^2n^{2\tau}
\bigr\} \leq \exp \bigl(-C_2n^{\min \{\eta\gamma,{(1-w)\gamma}/{\varrho_2},2\tau
,{\gamma(1-w)}/{6} \}} \bigr).
\]
\end{pn}

From Proposition~\ref{pn2}, we can find that the quantities on the right-hand
side are decreasing as $w$ is increasing. Thus, the optimal $w=\frac
{\kappa}{\varrho_1}$ is the same as the one for the best dimensionality
$p$ discussed previously.
Hence, the optimal bandwidth in our screening procedure is $w=\frac
{\kappa}{\varrho_1}$, which is quite sensible because intuitively the
smoother each $f_j$ is, the larger the bandwidth is allowed. The
corresponding upper bound for $\mathbb{P}\{\ell_j(0)\geq\half
c_1^2K_3^3n^{2\tau}\}$ is given in the following corollary.

\begin{cy}\label{cy2}
Under assumptions \textup{(A.1)}--\textup{(A.2)} and \textup{(A.4)}--\textup{(A.6)},
suppose $\max_{j\notin\mathcal{M}_*}\|f_j\|_\infty=O(n^{-\eta})$ for some $\eta>\frac
{(2\varrho_1+1)\kappa}{2\varrho_1}$. Pick $w=\frac{\kappa}{\varrho_1}$,
$\tau\in(\max\{\frac{1}{2}-\eta-\frac{\kappa}{2\varrho_1},0\},\frac
{1}{2}-\frac{(\varrho_1+1)\kappa}{\varrho_1})$ and $\xi>\frac{(\varrho
_1+2)\kappa}{\varrho_1}$. If $\inf_{u\in[a,b]}\mathbb{E}(Y^2|X_j=u)\geq
\rho$ for some positive $\rho$ holds for any $j\notin\mathcal{M}_*$,
then there exists a uniform positive constant $C_3$ such that for any
$j\notin\mathcal{M}_*$,
\[
\mathbb{P} \bigl\{\ell_j(0)\geq \halff c_1^2K_3^2n^{2\tau}
\bigr\}\leq p\exp \bigl(-C_3n^{\min \{\eta\gamma,{(\varrho_1-\kappa)\gamma
}/({\varrho_1\varrho_2}),2\tau,{(\varrho_1-\kappa)\gamma}/({6\varrho
_1}) \}} \bigr).
\]
\end{cy}

By noting that
\begin{eqnarray*}
|\widehat{\mathcal{M}}_{\gamma_n}|&=& \sum_{j\in\mathcal{M}_*}I
\biggl\{\ell _j(0)\geq {\half}c_1^2K_3^2n^{2\tau}
\biggr\}+\sum_{j\notin\mathcal
{M}_*}I \biggl\{\ell_j(0)\geq
\half c_1^2K_3^2n^{2\tau}
\biggr\}
\\
&\leq &  s+\sum_{j\notin\mathcal{M}_*}I \biggl\{
\ell_j(0)\geq \half c_1^2K_3^2n^{2\tau}
\biggr\},
\end{eqnarray*}
we have
$ \mathbb{P}(|\widehat{\mathcal{M}}_{\gamma_n}|>s) \leq \sum_{j\notin
\mathcal{M}_*}\mathbb{P}\{\ell_j(0)\geq \half c_1^2K_3^2n^{2\tau}\}$.
Hence, from Corollary~\ref{cy2}, we obtain the following theorem for the size
of $\widehat{\mathcal{M}}_{\gamma_n}$.

\begin{tm}
\label{tm2}
Under assumptions \textup{(A.1)}--\textup{(A.2)} and \textup{(A.4)}--\textup{(A.6)},
suppose $\max_{j\notin\mathcal{M}_*}\|f_j\|_\infty=O(n^{-\eta})$ for
some $\eta>\frac{(2\varrho_1+1)\kappa}{2\varrho_1}$. Pick $w=\frac
{\kappa}{\varrho_1}$, $\xi>\frac{(\varrho_1+2)\kappa}{\varrho_1}$ and
$\gamma_n=\frac{1}{2}c_1^2K_3^2n^{2\tau}$ for some $\tau\in(\max\{\frac
{1}{2}-\eta-\frac{\kappa}{2\varrho_1},0\},\frac{1}{2}-\frac{(\varrho
_1+1)\kappa}{\varrho_1})$. If $\inf_{u\in[a,b]}\mathbb{E}(Y^2|  X_j=u)\geq
\rho$ for some positive $\rho$ holds for any $j\notin\mathcal{M}_*$, then
\[
\mathbb{P} \bigl(|\widehat{\mathcal{M}}_{\gamma_n}|>s \bigr) \leq p\exp
\bigl(-C_3n^{\min \{\eta\gamma,{(\varrho_1-\kappa)\gamma}/({\varrho
_1\varrho_2}),2\tau,{(\varrho_1-\kappa)\gamma}/({6\varrho_1}) \}
} \bigr),
\]
where $C_3$ is given in Corollary~\ref{cy2}.
\end{tm}

This theorem shows that our screening procedure well controls
the set size of the recruited variables. With large probability,
the number of the recruited variables is not larger than the true size $s$.
From Theorems \ref{tm1} and \ref{tm2}
with $w=\frac{\kappa}{\varrho_1}$, we have that
$\mathbb{P}(\widehat{\mathcal{M}}_{\gamma_n}=\mathcal{M}_*)\rightarrow
1$ as $n\rightarrow\infty$ provided that
$\log p=o(n^{\min\{\eta\gamma,{(\varrho_1-\kappa)\gamma}/({\varrho
_1\varrho_2}),2\tau,{(\varrho_1-\kappa)\gamma}/({6\varrho_1}),1-2\kappa
-{\kappa}/{\varrho_1},({1}/{2}-\kappa-{\kappa}/{\varrho
_1}-\tau)\gamma\}})$.
This selection consistency property demonstrates that our approach
performs very well by distinguishing the true contributing variables
from false ones under condition~(\ref{eqfa}).
Aiming to obtain the optimal diverging rate for $p$, we can select
\[
\tau=\cases{
\displaystyle \frac{\gamma}{\gamma+2}\biggl(\frac{1}{2}-
\kappa-\frac{\kappa}{\varrho
_1}\biggr),&\quad $\displaystyle\mbox{if }\eta>\frac{1}{\gamma+2}+
\frac{\gamma\kappa}{\gamma
+2}+\frac{(\gamma-2)\kappa}{2(\gamma+2)\varrho_1}$;
\vspace*{3pt}\cr
\displaystyle\frac{1}{2}-\eta-\frac{\kappa}{2\varrho_1}+\varsigma, & \quad$\displaystyle\mbox{if }\eta \leq
\frac{1}{\gamma+2}+\frac{\gamma\kappa}{\gamma+2}+\frac{(\gamma
-2)\kappa}{2(\gamma+2)\varrho_1}$,
}
\]
where $\varsigma$ can be chosen to be positive and converging to $0$ as
$n\to \infty$. Hence, $\mathbb{P}(\widehat{\mathcal{M}}_{\gamma
_n}=\mathcal{M}_*)\rightarrow1$ as $n\rightarrow\infty$ provided that
\[
\log p=\cases{
o \bigl(n^{\min \{{(\varrho_1-\kappa)\gamma}/({\varrho_1\varrho
_2}),{\gamma} (1-2\kappa-{2\kappa}/{\varrho_1}
)/({\gamma+2}),{(\varrho_1-\kappa)\gamma}/({6\varrho_1}) \}} \bigr),\vspace*{3pt}\cr
\qquad\displaystyle\mbox{if }\eta>\frac{1+\gamma\kappa}{\gamma+2}+\frac{(\gamma-2)\kappa
}{2(\gamma+2)\varrho_1};
\vspace*{3pt}\cr
o \bigl(n^{\min \{{(\varrho_1-\kappa)\gamma}/({\varrho_1\varrho
_2}),{(\varrho_1-\kappa)\gamma}/({6\varrho_1}), (\eta-\kappa-
{\kappa}/({2\varrho_1}) )\gamma \}} \bigr),\vspace*{3pt}\cr\qquad\displaystyle\mbox{if }\eta\leq\frac
{1+\gamma\kappa}{\gamma+2}+
\frac{(\gamma-2)\kappa}{2(\gamma+2)\varrho
_1}.}
\]
More specifically,\vspace*{1pt} if the
response variable $Y$ has a compact support which means $\gamma=\infty
$, the above selection consistency holds if $\log p=o(n^{1-2\kappa-{2\kappa}/{\varrho_1}})$. Additionally, the smoothness of the
projections $f_j$'s also affects the allowable dimensionality. When all
$f_j\in C^{\infty}(\mathcal{X})$ implying that $\varrho_1=r=\infty$,
the allowable dimensionality turns out to be $\log p=o(n^{1-2\kappa})$.
If $Y$ follows a normal or sub-Gaussian distribution that $\gamma=2$
and $\eta=\infty$ which can be guaranteed by partial orthogonal
condition [\citet{HuangHorowitzWei2010AOS}], the selection consistency
holds if $\log p=o(n^{\min\{{1}/{2}-\kappa-{\kappa}/{\varrho
_1},{1}/{3}-{\kappa}/({3\varrho_1})\}})$.
It is worthwhile to note that though we show that our approach can
identify the set
of contributing variables with probability tending to 1, practical
performances can vary because first the results are valid
asymptotically, and second choosing the threshold level $\gamma_n$ to
achieve the perfect variable selection is difficult.

\section{Applications to some special models}\label{sec3}

Our local independence feature screening approach does not require a
specific form of the underlying model. Now we elaborate how the
proposed approach can be applied in three families of popular
nonparametric and semiparametric models: the nonparametric additive
models, the single-index models and
multiple-index models, and varying coefficient models; and we also
compare our results with existing ones.

\subsection{Nonparametric additive models}\label{sec31}
The
nonparametric additive model introduced by \citet{Stone1985} has the form
$
Y=\sum_{j=1}^ps_j(X_j)+\varepsilon$,
where\break $s_1(\cdot),\ldots, s_p(\cdot)$ are unknown functions with zero
mean and $\mathbb{E}(\varepsilon|\bX)=0$. It is a special case of model
(\ref{eqm}) with $m( \mathbf{X})=\sum_{j=1}^ps_j(X_j)$. For this model, the
true model can be defined as $ \mathcal{M}_*=\{1\leq j\leq p:\mathbb
{E}\{s_j^2(X_j)\}>0\}$.
Recall\vspace*{1pt} that $f_j(x)=\mathbb{E}(Y|X_j=x)$ following the earlier
discussion in Section~\ref{s2}. If we identify the true model by
%
\begin{equation}\label{eqid}
\min_{j\in\mathcal{M}_*}\|f_j\|_\infty\geq
c_1n^{-\kappa}
\end{equation}
for some nonnegative $\kappa$, our local independence feature screening
procedure proposed in Section~\ref{s2} can be directly applied here to surely
identify the contributing explanatory variables.

Let us carefully compare our approach with the one in \citet{FangFengSongJASA2011} that specifically targets at the feature screening problem
for nonparametric additive models. The screening procedure of \citet{FangFengSongJASA2011} includes four steps. First, each $f_j$ is expanded by
using the B-spline basis functions $\{\Psi_{j,k}\}_{k=1}^\infty$ and
use the truncated\vspace*{1pt} version $\tilde{f}_{nj}=\sum_{k=1}^{d_n}\beta
_{j,k}\Psi_{j,k}$ to approximate the projection $f_j$. Second, to
estimate the coefficients $\beta_{j,k}$'s and obtain the corresponding
estimation of each projection $f_j$. Third, to estimate each $\mathbb
{E}\{f_j^2(X_j)\}$ by $n^{-1}\sum_{i=1}^n\hat{f}_{nj}^2(X_{ij})$ where
$\hat{f}_{nj}$ is the estimation of $f_j$ obtained in the second step.
Fourth, to screen features via the corresponding magnitudes of these
estimates. To ensure the sure screening property, they assume that each
$f_j$ belongs to the class of functions whose $r$th derivative
satisfies the
Lipschitz continuity of order $\theta\in(0,1]$ and $d=r+\theta>0.5$,
where $r$ is a nonnegative integer, and identify the true model by the condition
%
\begin{equation}\label{eqiden}
\min_{j\in\mathcal{M}_*}\mathbb{E} \bigl\{f_j^2(X_j)
\bigr\}\geq Cd_nn^{-2\tilde{\kappa}}
\end{equation}
for some positive constant $C$. Here, $0<\tilde{\kappa}<\frac
{d}{2d+1}$ and $d_n$ is the truncation parameter used in approximating
each $f_j$ which satisfies $d_n\geq {C}n^{{2\tilde{\kappa}}/({2d+1})}$ for some positive constant ${C}$. By Theorem~1 of\vspace*{1pt} \citet{FangFengSongJASA2011}, the sure screening property holds for their procedure
if $\log p=o(n^{1-4\tilde{\kappa}}d_n^{-3})$. When $d_n=O(n^{{2\tilde{\kappa}}/({2d+1})})$, (\ref{eqiden}) implies\vspace*{1pt} $\min_{j\in
\mathcal{M}_*}\mathbb{E}\{f_j^2(X_j)\}\geq Cn^{-{4d\tilde{\kappa}}/({2d+1})}$ for some positive constant $C$. Actually, if the
density of each $X_j$ is uniformly bounded away from zero, these
conditions are sufficient for the identification condition of our
approach given in (\ref{eqid}) to hold\vspace*{1pt} with $\kappa=\frac{2d\tilde
{\kappa}}{2d+1}$. \citet{FangFengSongJASA2011} also assume that the error
$\varepsilon$ satisfies $\mathbb{E}\{\exp(B|\varepsilon|)|\bX\}\leq C$
for some positive constants $B$ and $C$
which implies there exist two positive constants $b_1$ and $b_2$ such
that $\mathbb{P}(|\varepsilon|\geq u)\leq b_1\exp(-b_2u)$ for any
$u>0$. See Lemma~2.2 in \citet{Petrov1995}. On the other hand, they also
assume that $\|m\|_\infty\leq \widetilde{B}$ for some positive
constant~$\widetilde{B}$. These two conditions together imply that $\gamma=1$ in
our assumption~(A.5). In this case,\vspace*{0.5pt} the sure screening property\vspace*{1.8pt} of our
approach given in Theorem~\ref{tm1}\vspace*{1pt} holds if $ \log p=o(n^{{1}/{2}-\kappa-{\kappa}/{\varrho_1}})=o(n^{{1}/{2}
-{2d\tilde{\kappa}}(1+1/\varrho_1)/(2d+1)})$. When
$d+\frac{1}{2}>(10+6d-\frac{2d}{\varrho_1})\tilde{\kappa}$, their
procedure can handle faster diverging $p$ than that of ours; otherwise,
our method is stronger than theirs. This shows that our procedure can
handle faster diverging $p$ when the signal strength level is weak,
that is, $\kappa$ is large.

In the specific case when the number of basis functions
$d_n=O(n^{{1}/({2d+1})})$ which leads to the optimal rate for the
B-spline approach [\citet{Stone1985}], their approach can handle the
dimensionality $ \log p=o(n^{1-4\tilde{\kappa}-{3}/({2d+1})})$.
For such a selection of $d_n$, the allowable dimensionality of\vspace*{1.5pt} our
approach under which the sure screening property holds is $ \log
p=o(n^{{1}/{2}+(1+1/\varrho_1)\{1/(4d+2)-\tilde{\kappa}\}})$. To make
the approach of \citet{FangFengSongJASA2011} work with a high-dimensional
setting for such a selection of $d_n$, the smooth parameter $d$ should
be larger than~$1$ implying the existence of the first derivation of
each $f_j$, which is not required in our approach.
From this point of view, our approach can handle the situation where
each nonparametric marginal projection $f_j$ does not have the first
derivative but being just continuous. When $d>1$, the parameter $r$ in
(A.1) and~(A.2) satisfies $1\leq r< d\leq r+1$. If the minimum signal
does not diminish to $0$ [\citet{LinZhang2006AOS,HuangHorowitzWei2010AOS}], then $\tilde{\kappa}=0$. In this case, our approach and
their approach can handle\vspace*{1.5pt} the nonpolynomial dimensionality $\log
p=o(n^{{1}/{2}+(1+1/r)/(4d+2)})$ and
$\log p=o(n^{1-{3}/({2d+1})})$, respectively. If $2d\leq6+{r}^{-1}$,
our approach allows faster diverging $p$ than that of their approach;
otherwise, their result is stronger than ours. This can be viewed as a
price paid for our approach by allowing weaker requirement on the
continuity of each $f_j$ and without requiring $m(\bX)$ to be bounded.

We note that the above diverging rates comparison is established in a
case in favor of the approach of \citet{FangFengSongJASA2011} by using\vspace*{1pt}
their identification condition with the smallest $d_n$. If $d_n$
diverges faster than $n^{{2\tilde{\kappa}}/({2d+1})}$, the
parameter $\kappa$ appeared in our identification condition will be
smaller than $\frac{2d\tilde{\kappa}}{2d+1}$ and the allowable
dimensionality of our approach will be improved. Additionally, our
approach has a very good control of the size of the recruited
variables. From Theorem~\ref{tm2}, the set of the recruited variables is not
larger than the true contributing covariates with large probability,
which together with Theorem~\ref{tm1}, imply the selection consistency of our
approach in nonparametric additive models.

\subsection{Single-index and multiple-index models}\label{sec32}
The single-index
model that is recognized as a particularly useful variation of the
linear regression model has the
form
$Y=s(\bbeta^{\mathrm{T}} \bX)+\varepsilon$,
where $s(\cdot)$ is the conditional mean function that is not
explicitly specified; see \citet{Brillinger1983} for more details. This kind
of models is a special case of model (\ref{eqm}) with $m(\bX)=s(\bbeta
^{\mathrm{T}} \bX)$.
Since single-index model \mbox{requires} identifiability condition
[\citet{Brillinger1983}], parameters in marginal single-index models is not
identifiable.
Therefore, marginal estimator-based ranking procedure cannot be applied
to handle the feature screening problem in the single-index model.

In this case, our local independence feature screening approach
conveniently applies,
because our marginal empirical likelihood based approach does not
require estimating parameters in the marginal models.
Since our marginal empirical likelihood approach is capable of
assessing $\mathbb{E}(Y|X_j)\equiv 0$, identifying the parameter in a
marginal single-index model is not necessary for our local independence
feature screening approach.
Specifically, if we identify the true model by $ \min_{j\in\mathcal
{M}_*}\|f_j\|_\infty\geq c_1n^{-\kappa}$, then the local independence
feature screening procedure and its properties discussed in
Section~\ref{s2}
also directly apply here. The effective application of our approach in
single-index models demonstrates that the marginal empirical likelihood
based approach is advantageous in independence feature screening. Such
a merit is due to the new insight of the marginal empirical likelihood
approach for screening variables by assessing the evidence against the
null hypothesis that the explanatory variable is not contributing marginally.

More generally, we may consider the screening for multiple-index model
$
Y=s(\bbeta_1^{\mathrm{T}}\bX,\ldots,\bbeta_K^{\mathrm{T}}\bX)+\varepsilon$,
where $K$ is a known integer, $\bbeta_k$ $(k=1,\ldots,K)$ are sparse
vectors of unknown parameters, and $s$ is an unknown function.
The marginal condition for the $j$th component of $\bX$ given by (\ref
{eqmj}) still applies in the multiple-index models.
Nevertheless, in multiple-index models, identification is also an
issue, which is actually more complicated and challenging than that in
the single-index models when considering the marginal contributions.
Since our concern in this model can also be transformed to assessing
$\mathbb{E}(Y|X_j)\equiv 0$ or not,
the local independence feature screening procedure given in the last
section can also be applied.

\subsection{Varying coefficient models}\label{sec33}
Varying coefficient model is useful for studying the variable-dependent
effects in the regressions. Many methods have been proposed for
estimation of this model. See, for example, \citet{FanZhang2000} for
the local polynomial smoothing method, \citet{HuangWuZhou2000} and \citet{QuLi2006} for basis expansion and spline method. The varying
coefficient model has the following form:
%
\begin{equation}\label{eqvarying}
Y=\bX^{\mathrm{T}} \bbeta(Z)+\varepsilon,
\end{equation}
where $\bX=(X_1,\ldots,X_p)^{\mathrm{T}}$ is a $p\times1$ vector of explanatory
variables, $Z$ is a scalar variable that takes values on a compact
interval $\mathcal{Z}$, and $\varepsilon$ is the error satisfies
$\mathbb{E}(\varepsilon|\bX,Z)=0$ almost surely. Here, $\bbeta(z)=(\beta
_1(z),\ldots,\beta_p(z))^{\mathrm{T}}$ is unknown but smooth in $z$. One way to
screen explanatory variables could be ignoring the impact due to $Z$,
and applying some marginal approaches. However, it is not difficult to
see that there are situations, for example, when $\int \beta(z)\,dz=0$,
that is, the parameter effect is zero in an average, univariate
screening procedure ignoring $Z$ will not be able to identify the
component in $\bX$ even if it is contributing.

To overcome such a difficulty, we may adjust the marginal nonparametric
regression for varying coefficient models to incorporate the effect of $Z$.
The marginal version of (\ref{eqvarying})
is
%
\begin{equation}
\label{eqmvary}
Y=a_j(Z)+X_jb_{j}(Z)+\tilde{
\varepsilon}_j,
\end{equation}
where $\mathbb{E}(\tilde{\varepsilon}_j|X_j,Z)=0$ almost surely for
$j=1,\ldots,p$. Here, $b_{j}(Z)$ can be interpreted as the marginal
contribution of $X_j$ in explaining $Y$ via $Z$. It can be shown that
$b_j(Z)=\frac{\operatorname{cov}(X_j,Y|Z)}{\operatorname{var}(X_j|Z)}$.
Thus, the marginal effect $b_j(z)=0$ for any $z$ is equivalent to
$\operatorname{cov}(X_j,Y|Z=z)=0$ for any $z$.
In the simple case when $\mathbb{E}(Y|Z)=0$ or some constant free of
$Z$, it is equivalent to $\operatorname{cov}(X_j,Y|Z)=\mathbb
{E}(X_jY|Z)\equiv 0$, which essentially share the same form $\mathbb
{E}(Y|X_j)\equiv 0$ as in the local independence featuring screening.
In the general situation, $\mathbb{E}(Y|Z)\neq 0$ so that one needs to
assess $0\equiv \operatorname{cov}(X_j,Y|Z)=\mathbb{E}[X_j \{Y-\mathbb
{E}(Y|Z)\}|Z]$ with the nuisance function $\mathbb{E}(Y|Z)$ estimated.
For a kernel function $\widetilde{\mathcal{K}}(\cdot)$, $\mathbb
{E}(Y|Z=z)$ can be estimated by
%
\begin{equation}
\label{eqnw}
\widehat{\mathbb{E}}(Y|Z=z)=\frac{n^{-1}\sum_{i=1}^n\widetilde{\mathcal
{K}}_{\tilde{h}}(Z_i-z)Y_i}{n^{-1}\sum_{i=1}^n\widetilde{\mathcal
{K}}_{\tilde{h}}(Z_i-z)},
\end{equation}
where $\widetilde{\mathcal{K}}_{\tilde{h}}(u)=\tilde{h}^{-1}\widetilde
{\mathcal{K}}({u}\tilde{h}^{-1})$ with bandwidth $\tilde{h}$. Let
$\widetilde{Y}_i=Y_i-\widehat{\mathbb{E}}(Y|Z=Z_i)$, then the marginal
empirical likelihood is constructed as
%
\begin{eqnarray}\label{eqmalv}
&&\quad\mathrm{EL}_j(z,0)
=\sup \Biggl\{\prod_{i=1}^nw_i:w_i
\geq0,\sum_{i=1}^nw_i=1,\sum
_{i=1}^nw_i
\mathcal{K}_h(Z_{i}-z)X_{ij}\widetilde
Y_i=0 \Biggr\}.\hspace*{-12pt}
\end{eqnarray}
We then propose using
%
\begin{equation}
\label{eqelv}
\ell_j(0)=\sup_{z\in\mathcal{Z}_n}
\ell_j(z,0),
\end{equation}
where $\ell_j(z,0)=-2\log\{\mathrm{EL}_j(z,0)\}-2n\log n$, and $\mathcal
{Z}_n$ is a partition of $\mathcal{Z}$. In practice, a natural choice
is to evaluate the statistic (\ref{eqelv}) by $\ell_j(0)=\max_{1\leq
i\leq n}\ell_j(Z_{i},0)$ where $\{Z_i\}_{i=1}^n$ are the $n$
observations of variable $Z$. This new $\ell_j(0)$ can be used for
local independence feature screening for varying coefficient models
analogously to that in Section~\ref{s2}.

The analogous assumptions corresponding to (A.1)--(A.5) in Section~\ref{s2} are
given as follows.
\begin{longlist}[(A.5)$'$]
\item[(A.1)$'$]  For each $j=1,\ldots,p$, let $\tilde{f}_j(z)=\textrm
{cov}(X_j,Y|Z=z)$. Assume $\{\tilde{f}_j\}_{j=1}^p$ belong to
$C^r(\mathcal{Z})$. If $r=0$, we assume that $\tilde{f}_{j}$'s satisfy
the Lipschitz condition with order $\alpha\in(0,1]$, that is, $ |\tilde
{f}_{j}(s)-\tilde{f}_{j}(t)|\leq \widetilde{K}_1|s-t|^\alpha$ for any $s,t\in\mathcal{Z}$, where $\widetilde{K}_1$ is a positive
constant uniformly for any $j=1,\ldots,p$. In addition, there exists a
constant $\widetilde{K}_2$ such that $|\tilde{f}_{j}^{(r)}(z)|\leq
\widetilde{K}_2$ for any $z\in\mathcal{Z}$ and $j=1,\ldots,p$.

\item[(A.2)$'$]  The density function $g$ of $Z$ satisfies $0<\widetilde
{K}_3\leq g(z)\leq \widetilde{K}_4<\infty$ on $\mathcal{Z}$. In
addition, we assume that $g$ belongs to $C^r(\mathcal{Z})$ for the $r$
given in (A.1)$'$ and $|g^{(r)}(z)|\leq \widetilde{K}_5$ for any $z\in
\mathcal{Z}$.

\item[(A.3)$'$]  There exist nonnegative constants $\tilde{c}_1>0$ and $\kappa\in
[0,\frac{\max\{r,\alpha\}}{2\max\{r,\alpha\}+2})$ such that $ \min_{j\in
\mathcal{M}_*}\|\tilde{f}_j\|_\infty\geq \tilde{c}_1n^{-\kappa}$, where
$r$ and $\alpha$ are specified in (A.1)$'$.

\item[(A.4)$'$] There exists some positive constant $\xi$ such that $\|\mathcal
{Z}_n\|=n^{-\xi}$.

\item[(A.5)$'$]  There exist positive constants $\widetilde{K}_6, \widetilde
{K}_7, \gamma_1$ and $\gamma_2$ such that $
\mathbb{P}(|Y|\geq u)\leq \widetilde{K}_6\exp(-\widetilde{K}_7u^{\gamma_1})$
and
$\mathbb{P}(|X_j|\geq u)\leq \widetilde{K}_6\exp(-\widetilde
{K}_7u^{\gamma_2}) $ for any $u>0$ and $j=1,\ldots,p$.

Assumption (A.5)$'$ and Lemma~2 in \citet{ChangTangWu2013b} yield that
$ \mathbb{P}(|X_jY|\geq u)\leq 2\widetilde{K}_6\exp(-\widetilde
{K}_7u^\gamma) $
for any $u>0$ and $j=1,\ldots,p$, where $\gamma=\frac{\gamma_1\gamma
_2}{\gamma_1+\gamma_2}$.
To investigate\vspace*{1.5pt} the theoretical properties of (\ref{eqelv}) with
nonparametric estimation (\ref{eqnw}), we need the following two extra
conditions:

\item[(A.7)]  $\mathbb{E}(Y|Z=z)$ belongs to $C^r(\mathcal{Z})$, where $r$ is
given in (A.1)$'$. If $r=0$, we assume $\mathbb{E}(Y|Z=z)$ satisfies the
Lipschitz condition with order $\alpha$ where $\alpha$ is specified in (A.1)$'$.

\item[(A.8)]  For $r$ specified in (A.1)$'$, if $r\geq1$, the kernel function
$\widetilde{\mathcal{K}}(\cdot)$ is of order $r$. If $r=0$, the kernel
function satisfies $\widetilde{\mathcal{K}}(u)\geq0$ for any $u$ and
$\int\widetilde{\mathcal{K}}(u)\,du=1$.
\end{longlist}
For $\gamma=\frac{\gamma_1\gamma_2}{\gamma_1+\gamma_2}$, we define
$\varrho_1$, $\varrho_2$ and $\delta$ as (\ref{eqparameters}). In
addition, let $\delta_1=\max\{\frac{2}{\gamma_1}-1,0\}$. Meanwhile, we
assume that the bandwidths $h$ and $\tilde{h}$ in constructing marginal
empirical likelihood (\ref{eqmalv}) and the NW estimator (\ref{eqnw})
for $\mathbb{E}(Y|Z=z)$ satisfy $ h\asymp n^{-w}$ and $\tilde
{h}\asymp n^{-\phi}$, where $w$ and $\phi$ will be specified later. The
property of the marginal empirical likelihood for varying coefficient
models are given in the following theorem.

\begin{tm}\label{cy3}
Under assumptions \textup{(A.1)}$'$--\textup{(A.5)}$'$, \textup{(A.6)} and
\textup{(A.7)}--\textup{(A.8)}, pick $w\in[\frac{\kappa}{\varrho_1},\frac{1}{2}-\kappa
)$, $\phi\in(\frac{\kappa}{\varrho_1},1-2\kappa)$, $\xi>\kappa+2w$, and
$\gamma_n=\frac{1}{2}\tilde{c}_1^2\widetilde{K}_3^2n^{2\tau}$ for some
$\tau\in(0,\frac{1}{2}-\kappa-w)$, then there exists a uniform positive
constant $C_4$ such that
%
%
\begin{eqnarray*}
 \mathbb{P} (\mathcal{M}_*\subset\widehat{\mathcal{M}}_{\gamma
_n} )
&\geq &
1 -s\exp \bigl\{-C_4n^{({1}/{2}-\kappa-w-\tau)\gamma} \bigr\}-s\exp \bigl
\{-C_4n^{(\phi\varrho_1-\kappa)\gamma_2}\bigr\}
\\
&&{}-s\exp \bigl(-C_4n^{\min \{1-2\kappa-w,\frac{1-\kappa-w}{1+\delta}
\}} \bigr)
\\
&&{}-s\exp \bigl(-C_4n^{\min \{\frac{\{3-2(\phi+\kappa+w+\tau)\}\gamma
_2}{(2+2\delta_1)\gamma_2+2},\frac{\{2-\phi-2(\kappa+w+\tau)\}\gamma
_2}{\gamma_2+2} \}} \bigr)
\\
&&{}-s\exp \bigl(-C_4n^{\min \{\frac{(1-2\kappa-\phi)\gamma_2}{\gamma
_2+2},\frac{(1-\kappa-\phi)\gamma_2}{(1+\delta_1)\gamma_2+1} \}} \bigr).
\end{eqnarray*}
%
\end{tm}

Theorem~\ref{cy3} provides a general result for the sure screening
property of the marginal empirical likelihood for varying coefficient
models. The first two terms and the fourth term on the right-hand side
of above inequality are the same as those in Theorem~\ref{tm1}. The extra terms
are due to the kernel estimation of $\mathbb{E}(Y|Z)$. The analogues of
Proposition~\ref{pn1} and Theorem~\ref{tm2} are also valid for varying coefficient
models using the above marginal empirical likelihood.

\citet{Fanet2014}, \citet{Liuetal2014} and
\citet{Songetal2014} also consider the feature screening for ultra-high dimensional
varying coefficient models.
\citet{Fanet2014} and \citet{Liuetal2014} considered the same
model as (\ref{eqvarying}) while \citet{Songetal2014} allows both
$Y$ and $\bX$ to depend on $Z$.
\citet{Fanet2014} estimate $a_j(\cdot)$ and $b_j(\cdot)$ in (\ref
{eqmvary}) simultaneously via
the B-spline basis functions expansion approach.
\citet{Fanet2014} propose an estimator $\hat{u}_j$ for $\mathbb
{E}\{\frac{\operatorname{cov}^2(X_j,Y|Z)}{\operatorname{var}(X_j|Z)}\}$ for each
$j=1,\ldots,p$ measuring the marginal contribution of $X_j$.
Then they propose to select the covariates via ranking $|\hat{u}_j|$'s.
\citet{Liuetal2014}
study conditional Pearson correlation as a measure of marginal
contribution for varying coefficient models. For each $j=1,\ldots,p$
and $z$, they estimate conditional Pearson\vspace*{1.5pt} correlation $\rho
(X_j,Y|Z=z)$ via kernel smoothing method and construct an estimation
$\hat{\rho}_j^*$ for $\mathbb{E}\{\rho^2(X_j,Y|Z)\}$. Then they use the
magnitude of $|\hat{\rho}_j^*|$ to determine whether $X_j$ is an
important explanatory variable or not.
The main idea of \citet{Songetal2014} is similar to that of \citet{Fanet2014}. In the sequel, we carefully compare our procedure
with those proposed in \citet{Fanet2014} and \citet{Liuetal2014}.

\citet{Fanet2014} identify $\mathcal{M}_*$ via $ \min_{j\in\mathcal
{M}_*}\mathbb{E}\{\operatorname{cov}^2(X_j,Y|Z)\}\geq Cd_nn^{-2\widetilde
{\kappa}} $ for some positive constants $C$ and $\tilde{\kappa
}<\frac{2d+1}{8d+10}$. Here, $d_n$ is the number of approximation terms
to $a_j(\cdot)$ and $b_j(\cdot)$ in the estimation step. Details of
$d_n$ can be found in Section~\ref{sec31} for discussion of additive models. To
guarantee the validity of the approach in \citet{Fanet2014},
$d_n\geq {C}n^{{2\tilde{\kappa}}/({2d+1})}$ is required for some
positive constant ${C}$ and $d=r+\theta$, where $r$ is an integer and
$\theta\in(0,1]$ are employed to describe the smoothness of each $a_j$
and $b_j$, that is, the $r$th derivatives of all $a_j$ and $b_j$ are
Lipschitz continuous of order $\theta$. In \citet{Liuetal2014}, the
identification condition is given by $ \min_{j\in\mathcal{M}_*}\mathbb
{E}\{\rho^2(X_j,Y|Z)\}\geq Cn^{-2\bar{\kappa}} $ for some positive
constants $C$ and
$\bar{\kappa}<\hbar$. Here, $\hbar<\frac{1}{3}$ is a parameter employed
to describe the decay rate of the bandwidth for the kernel smoothing
step in their procedure, that is, the bandwidth is chosen as
$O(n^{-\hbar})$. Based on the moments condition and the assumptions
$\inf_{z\in\mathcal{Z}}\min_{1\leq j\leq p}\operatorname{var}(X_j|Z=z)>0$ and
$\inf_{z\in\mathcal{Z}}\operatorname{var}(Y|Z=z)>0$, the identification
condition in \citet{Liuetal2014} is essentially equivalent to $ \min_{j\in\mathcal{M}_*}\mathbb{E}\{\operatorname{cov}^2(X_j,Y|Z)\}\geq Cn^{-2\bar
{\kappa}} $ for some positive constant $C$.

From three aspects, we compare our identification condition and
theoretical results with those of \citet{Fanet2014} and \citet{Liuetal2014}. First, if the density of $Z$ is uniformly bounded away
from zero and infinity on its support $\mathcal{Z}$, their $L_2$-type
requirement is a sufficient condition for ours proposed in assumption
(A.3)$'$. But their identification conditions rule out the case where
$\operatorname{cov}(X_j,Y|Z=z)$ only contribute largely at several local
small intervals on $\mathcal{Z}$. Second, our method works for weaker
signal strength than theirs. \citet{Fanet2014} can handle the case
that $[\mathbb{E}\{\operatorname{cov}^2(X_j,Y|Z)\}]^{1/2}\geq Cn^{-{d\tilde{\kappa}}/({2d+1})}$ for some $\tilde{\kappa}<\frac
{2d+1}{8d+10}$. Therefore, the weakest signal strength can be handled
by their method cannot decay at a rate faster than $O(n^{-{d}/({8d+10})})$. On the other hand, the weakest signal strength our
method can handle is at the rate of
$O(n^{-{\varrho_1}/({2\varrho_1+2})})$. By noting that $r<d\leq r+1$,
we have $\frac{\varrho_1}{2\varrho_1+2}>\frac{d}{8d+10}$ which implies
our method can accommodate weaker signal strength than \citet{Fanet2014}. If condition (C4) of \citet{Liuetal2014} holds, the parameter
$r$ in our notation system is not smaller than $2$. Thus, $\frac{\varrho
_1}{2\varrho_1+2}\geq\frac{1}{3}$ which implies our method can also
accommodate weaker signal strength than that of \citet{Liuetal2014}.
Third, we compare the nonpolynomial dimensionality allowed by different
methods. As our theoretical assumptions are close to that proposed in
\citet{Liuetal2014}, we compare our result with theirs. When $\gamma
_1=\gamma_2=r=2$, which are assumed in \citet{Liuetal2014}, from
Corollary~3, our method can handle $ \log p=o(n^{\min\{4\phi-2\kappa
,{1}/{2}-\kappa-{\phi}/{2}\}}) $
when $\tau$ is chosen to be close enough to zero and $w=\frac{\kappa
}{\varrho_1}=\frac{\kappa}{2}$, where the result for the method of \citet{Liuetal2014} is $\log p=o(n^{\hbar-\kappa})$ for some $\hbar<\frac
{1}{3}$. Notice that $\kappa<\hbar<\frac{1}{3}$ in their setting. If we
choose $\phi=\frac{1}{5}$, then $\min(4\phi-2\kappa,\frac{1}{2}-\kappa
-\frac{\phi}{2})=\frac{2}{5}-\kappa>\hbar-\kappa$, which means our
procedure can accommodate faster diverging $p$ than theirs.

\section{Iterative feature screening}\label{sec4}

It is possible that some predictors are margi\-nally unrelated but
jointly related to the response as illustrated by Example~4.2.2 of \citet{FanLv2008JRSSB}. As observed in the literature, marginal utility-based
feature screening methods may miss this type of predictors badly [\citet{FanLv2008JRSSB,FSW2009}]. To overcome this difficulty,
several versions of iterative feature screening have been proposed. \citet{FanLv2008JRSSB} proposed to regress the response on the recruited
predictors and use the regression residual as the new ``response'' to
recruit further from the remaining predictors. While recruiting
additional predictors, \citet{FSW2009} consider the joint
model with both the recruited predictors and each additional predictor
and use the conditional contribution of each additional predictor given
those predictors that are already selected to recruit further from the
remaining predictors. Both versions need to fit a joint model on the
recruited predictors
with or without an additional feature where a parametric model
specification is required. However our proposed empirical
likelihood-based screening is for a general nonparametric model (\ref
{eqm}) and, in this sense, model-free. Consequently, the aforementioned
two versions of iterative feature screening cannot be extended to our
case. Next, we will borrow the idea of \citet{Zhuetal2011JASA} and propose an
iterative version for our proposed empirical likelihood-based screening.

Next, we detail our iterative screening procedure.
\begin{longlist}[\textit{Step} 4.]
\item[\textit{Step} 1.] We first apply our proposed empirical likelihood based
screening to $\{(\bX_i, Y_i), i=1, \ldots, n\}$ and denote the selected
set of predictors by $\widehat{\mathcal{A}}_1$.

\item[\textit{Step} 2.] Apply COSSO [\citet{LinZhang2006AOS}] to data $\{(\bX_{i\widehat{\mathcal{A}}_1}, Y_i), i=1,\ldots, n\}$
and use $\widehat{\mathcal{M}}_1$ to denote the subset of $\widehat{\mathcal{A}}_1$ that are retained by the COSSO.

\item[\textit{Step} 3.] For any $j\notin\widehat{\mathcal{M}}_1$, we regress $X_j$ on
the predictors with indices in $\widehat{\mathcal{M}}_1$ based on the
data $\{(\bX_{i\widehat{\mathcal{M}}_1}, X_{ij}), i=1, \ldots, n\}$ and
denote the residual by $\hat\varepsilon_{ij}$, $i=1, \ldots, n$. For
any $j\notin\widehat{\mathcal{M}}_1$, treat $\hat\varepsilon_{ij}$,
$i=1, \ldots, n$ as the pseudo predictor and apply our proposed
empirical likelihood-based screening to recruit a subset $\widehat
{\mathcal{A}}_2$ of predictors.

\item[\textit{Step} 4.] Apply COSSO to data $\{(\bX_{i\widehat{\mathcal{M}}_1\cup
\widehat{\mathcal{A}}_2}, Y_i), i=1,\ldots, n\}$ and use $\widehat
{\mathcal{M}}_2$ to denote the subset of $\widehat{\mathcal{M}}_1\cup
\widehat{\mathcal{A}}_2$ that are retained by the COSSO.

\item[\textit{Step} 5.] Repeat Steps 3 and 4 until $\widehat{\mathcal{M}}_k=\widehat
{\mathcal{M}}_{k-1}$ or the number of predictors in $\widehat{\mathcal
{M}}_k$ reaches some prespecified number.
\end{longlist}

\section{Numerical examples}\label{sec5}

We next demonstrate the performance of the local independence feature
screening methods using five examples with comparisons to appropriate
existing alternative approaches. In what follows, we denote our method
by EL when reporting the results. When implementing the local
independence feature screening approach, we select the bandwidth $h$ by
the cross-validation method [\citet{fan1996}], and the
Epanechnikov kernel $K(u)=\frac{3}{4}(1-u^2) I(|u|\leq 1)$ is applied
for the marginal regressions.

\begin{exa}\label{exa1}
This example is taken from Example~3 of \citet{FangFengSongJASA2011}. Data are generated from model
$Y=5g_1(X_1)+3g_2(X_2)+4g_3(X_3)+6g_4(X_4)+\sigma\varepsilon$ with
$g_1(x)=x$, $g_2(x)=(2x-1)^2$, $g_3(x)=\frac{\sin(2\pi x)}{2-\sin(2\pi
x)}$, and $g_4(x)=0.1\sin(2\pi x)+0.2\cos(2\pi x)+0.3\sin^2(2\pi
x)+0.4\cos^3(2\pi x)+\break 0.5 \sin^3(2\pi x)$. Here, predictors $X_j$'s are
i.i.d. random variables of $\operatorname{Unif}(0,1)$ distribution, and
$\varepsilon\sim N(0,1)$ is independent of $X_j$'s. We set $p=1000$.
Data sets of size $n=400$ are used. We apply the proposed screening
method to reduce the number of predictors from 1000 to 20. We consider
different signal noise ratios by varying $\sigma^2$ at four different
%
\begin{table}
\caption{Simulation results for Example~\protect\ref{exa1}}\label{tbex1}
\begin{tabular*}{\textwidth}{@{\extracolsep{\fill}}lccccc@{}}
\hline
$\bolds{\sigma^2}$ & \textbf{Method} & $\bolds{X_1}$ & $\bolds{X_2}$ & $\bolds{X_3}$ & $\bolds{X_4}$\\
\hline
1 &EL &100&  97 & 100 & 100\\
&\citet{FangFengSongJASA2011} &100&  90 & 100 & 100\\
 &\citet{Zhuetal2011JASA}&100& \phantom{0}1 & 100 & 100\\
&\citet{Lietal2012JASA} &100&27 & 100 & 100\\
& \citeauthor{ChangTangWu2013} (\citeyear{ChangTangWu2013}) &100& \phantom{0}1 & 100 & 100\\[3pt]
1.74 &EL &100&  96 & 100 & 100\\
&\citet{FangFengSongJASA2011} &100&  88 & 100 & 100\\
 &\citet{Zhuetal2011JASA}&100& \phantom{0}3 & 100 & 100\\
&\citet{Lietal2012JASA} &100& 24 & 100 & 100\\
&\citeauthor{ChangTangWu2013} (\citeyear{ChangTangWu2013}) &100& \phantom{0}1 & 100 & 100\\[3pt]
2&EL &100&  95 & 100 & 100\\
&\citet{FangFengSongJASA2011} &100&  87 & 100 & 100\\
 &\citet{Zhuetal2011JASA}&100& \phantom{0}2 & 100 & 100\\
&\citet{Lietal2012JASA} &100& 24 & 100 & 100\\
&\citeauthor{ChangTangWu2013} (\citeyear{ChangTangWu2013}) &100& \phantom{0}1 & 100 & 100\\[3pt]
3&EL &100&  93 & 100 & 100\\
&\citet{FangFengSongJASA2011} &100& 85 & 100 & 100\\
 &\citet{Zhuetal2011JASA}&100& \phantom{0}1 & 100 & 100\\
&\citet{Lietal2012JASA} &100& 20 & 100 & 100\\
&\citeauthor{ChangTangWu2013} (\citeyear{ChangTangWu2013}) &100& \phantom{0}1 & 100 & 100\\
\hline
\end{tabular*}\vspace*{3pt}
\end{table}
levels while \citet{FangFengSongJASA2011} chose a specific value of $\sigma
^2=1.74$. Table~\ref{tbex1} reports the frequency of the important
predictors being selected among 100 repetitions for different values of
$\sigma^2$. A comparison is made with \citet{FangFengSongJASA2011}. It is observed that the proposed empirical
likelihood-based local independence feature screening
performs similarly as the method proposed by \citet{FangFengSongJASA2011}. Contributing features $X_1$, $X_2$ and $X_4$ are selected by
both methods for all repetitions. Yet for feature $X_2$ with a slightly
weaker contribution, our method does slightly better.
\end{exa}

\begin{exa}\label{exa2}
This is a nonlinear predictor effect example with
heterogeneous conditional variance. Data are generated from model
$Y=-3h_1(X_1)+2.5 h_2(X_2)-2h_3(X_3)+1.5 h_4(X_4)+\sigma\varepsilon$
with\vspace*{1pt} $h_1(x)=(2x-1)^2$, $h_2(x)=\break \frac{\cos(2\pi x)}{2+\sin(2\pi x)}$,
$h_3(x)=\frac{\cos(2\pi x)}{2-\cos(2\pi x)}$, and\vspace*{1.5pt} $h_4(x)=\cos\{
(2x-1)\pi\}$. Here, $X_j$'s are independent and uniform over $[0,1]$,
$\varepsilon$ is independent of $X_j$'s and has normal distribution
with mean zero and its heterogeneous conditional variance is generated
by $\operatorname{var}(\varepsilon)=\frac{4}{x_1^2+x_2^2+x_3^2 +x_4^2}$. The
noise\vspace*{1pt} level is govern by $\sigma$ with different values in the
simulations. We set $p=1000$ and $n=100$. We apply the proposed
\begin{table}
\tabcolsep=0pt
\caption{Simulation results for Example~\protect\ref{exa2} with $\sigma^2$ controlling
the overall noise level}\label{tbexn}
\begin{tabular*}{\textwidth}{@{\extracolsep{4in minus 4in}}lccccccccc@{}}
\hline
&& \multicolumn{4}{c}{\textbf{Homogeneous}}& \multicolumn
{4}{c@{}}{\textbf{Heterogeneous}}\\
&& \multicolumn{4}{c}{\textbf{variance}}&\multicolumn
{4}{c@{}}{\textbf{variance}}\\[-4pt]
&& \multicolumn{4}{l}{\hrulefill} & \multicolumn{4}{l@{}}{\hrulefill}\\
$\bolds{\sigma^2}$&\textbf{Method} & $\bolds{X_1}$ & $\bolds{X_2}$ & $\bolds{X_3}$ & $\bolds{X_4}$&$\bolds{X_1}$ & $\bolds{X_2}$ & $\bolds{X_3}$
& $\bolds{X_4}$\\
\hline
0.5  &EL & 99&97&97&100& 83&90&81&94 \\
&\citet{FangFengSongJASA2011} &94&99&90&100& 76&86&78&92 \\
&\citet{Zhuetal2011JASA} & \phantom{0}3&\phantom{0}4&\phantom{0}4&\phantom{00}4&\phantom{0}4&\phantom{0}6&\phantom{0}2&\phantom{0}2\\
&\citet{Lietal2012JASA} & 36&59&38&\phantom{0}68&23&43&22&47\\
& \citeauthor{ChangTangWu2013} (\citeyear{ChangTangWu2013}) &\phantom{0}3&\phantom{0}4&\phantom{0}1&\phantom{00}2&\phantom{0}4&\phantom{0}6&\phantom{0}0&\phantom{0}0\\[3pt]
1.0 &EL &94&98&89&\phantom{0}99& 66&74&67&85 \\
&\citet{FangFengSongJASA2011} &88&95&86&\phantom{0}97& 53&67&58&84 \\
 &\citet{Zhuetal2011JASA} & \phantom{0}3&\phantom{0}4&\phantom{0}2&\phantom{00}3&\phantom{0}4&\phantom{0}4&\phantom{0}0&\phantom{0}4 \\
&\citet{Lietal2012JASA} &27&50&32&\phantom{0}60& 16&34&16&37\\
&\citeauthor{ChangTangWu2013} (\citeyear{ChangTangWu2013}) &\phantom{0}3&\phantom{0}5&\phantom{0}1&\phantom{00}2& \phantom{0}3&\phantom{0}5&\phantom{0}0&\phantom{0}1 \\[3pt]
1.5 &EL &88&95&87&\phantom{0}97& 57&61&48&73 \\
&\citet{FangFengSongJASA2011} & 81&92&81&\phantom{0}96&43&55&41&76\\
 &\citet{Zhuetal2011JASA} & \phantom{0}3&\phantom{0}6&\phantom{0}2&\phantom{00}2&\phantom{0}4&\phantom{0}4&\phantom{0}1&\phantom{0}2\\
&\citet{Lietal2012JASA} &22&44&26&\phantom{0}51& 14&24&13&26\\
&\citeauthor{ChangTangWu2013} (\citeyear{ChangTangWu2013}) & \phantom{0}3&\phantom{0}5&\phantom{0}0&\phantom{00}2&\phantom{0}3&\phantom{0}5&\phantom{0}0&\phantom{0}1 \\[3pt]
2.0 &EL & 84&86&82&\phantom{0}94&46&52&38&65 \\
&\citet{FangFengSongJASA2011} & 77&87&79&\phantom{0}93&34&50&35&61\\
 &\citet{Zhuetal2011JASA} &\phantom{0}4&\phantom{0}5&\phantom{0}0&\phantom{00}2& \phantom{0}4&\phantom{0}5&\phantom{0}1&\phantom{0}2\\
&\citet{Lietal2012JASA} &21&39&21&\phantom{0}47& 14&20&11&20\\
&\citeauthor{ChangTangWu2013} (\citeyear{ChangTangWu2013}) &\phantom{0}3&\phantom{0}5&\phantom{0}0&\phantom{00}1&\phantom{0}5&\phantom{0}4&\phantom{0}0&\phantom{0}2 \\
\hline
\end{tabular*}\vspace*{3pt}
\end{table}
screening method to reduce the number of predictors from 1000 to 20.
For comparison purposes, we also apply the methods on data generated
from a model with homogenous conditional variance while other settings
are the same. The results are reported in Table~\ref{tbexn} for the two
cases. There are a few observations. First, those methods incorporating
less local impact perform poorly in this nonlinear effect example,
demonstrating the importance and
substantial effect of the local feature of the marginal contributions
to the response variable. Second, in this example with heterogeneous
conditional variance, our method outperforms others, and is better than
the one of \citet{FangFengSongJASA2011}, especially when the  noise
level is relatively higher implying the signal is relatively weaker.
This is consistent with our theory and the finding from Example~\ref{exa1}, and
it shows that our method is advantageous for detecting nonlinear
effects. It also demonstrates that when the signal is weak, and when
the situation is more difficult due to high level and more complex
variations, our method delivers more promising results thanks to the
feature of the marginal empirical likelihood approach.\looseness=1
\end{exa}

\begin{exa}\label{exa3}
Data are generated from model $ Y=\beta_1X_1+\beta
_2X_2+\beta_3X_3+\beta_4X_4+\varepsilon$ with independent error
$\varepsilon\sim N(0,1)$. Jointly Gaussian covariates satisfy $\mathbb
{E}(X_j)=0$ and $\operatorname{var}(X_j)=1$ for $j=1,\ldots, p$ with $\operatorname{cov}(X_j,X_4)=\frac{1}{\sqrt{2}}$ for $j\ne 4$ and
$\operatorname{cov}(X_j,
X_{j'})=\frac{1}{2}$ if $j$ and $j'$ are distinct elements of $\{
1,\ldots,p\}\setminus \{4\}$. True regression coefficients are given by
$\beta_1=\beta_2=\beta_3=2$, $\beta_4=-3\sqrt{2}$, and $\beta_j=0$ for
$j>4$ such that $X_4$ is marginally independent of the response $Y$.
Yet $X_4$ is the most important predictor variable in the joint model.
This example is to illustrate that the iterative version of the
proposed screening procedure works effectively. We borrow the idea of
\begin{table}
\caption{Simulation results for Example~\protect\ref{exa3}}\label{tbex2}
\begin{tabular*}{\textwidth}{@{\extracolsep{\fill}}lcccccc@{}}
\hline
&&&&&& \textbf{Average} \textbf{for}\\
$\bolds{(n,p)}$& \textbf{Iterative} \textbf{screening} \textbf{method} & $\bolds{X_1}$ & $\bolds{X_2}$ & $\bolds{X_3}$ & $\bolds{X_4}$ &
  $\bolds{x_j}$\textbf{,} $\bolds{j\geq 5}$\\
\hline
(300, 1000) &EL & 100 & 100 & 100 & 100& 0.0351\\
&\citet{FangFengSongJASA2011} & 100 & 100 & 100 & 100& 0.1175\\
&\citet{Zhuetal2011JASA}&100 & 100 & 100 & 100&N/A\\
&\citet{Lietal2012JASA} &100 & 100 & 100 & 100&N/A\\
& \citeauthor{ChangTangWu2013} (\citeyear{ChangTangWu2013}) &100 & 100 & 100 & \phantom{0}99&0.0281\\[3pt]
(400, 2000) &EL& 100 & 100 & 100 & 100& 0.0141\\
&\citet{FangFengSongJASA2011} & 100 & 100 & 100 & 100& 0.0612\\
&\citet{Zhuetal2011JASA}&100 & 100 & 100 & 100&N/A\\
&\citet{Lietal2012JASA} &100 & 100 & 100 & 100&N/A\\
&\citeauthor{ChangTangWu2013} (\citeyear{ChangTangWu2013}) &100 & 100 & 100 & 100&0.0220\\
\hline
\end{tabular*}\vspace*{3pt}
\end{table}
\citet{Zhuetal2011JASA} to define the iterative version of our screening
procedure as laid out in Section~\ref{sec4}. \citet{FangFengSongJASA2011} only
considered the case with $(n=400,
p=1000)$ while we consider two cases: $(n=300, p=1000)$ and $(n=400,
p=2000)$. Simulation results over 100 repetitions are reported in
Table~\ref{tbex2}, where we report the frequency of important
predictors being selected and the average frequency of unimportant
predictors being selected. It shows that the iterative screening based
on empirical likelihood performs similarly as the nonparametric
screening proposed in \citet{FangFengSongJASA2011}. In terms of average
frequency of unimportant predictors being selected, our new method has
slight advantage.
\end{exa}

\begin{exa}\label{exa4}
In this example, we consider a single-index type
model. Data are generate from $Y=m(\bX)+\sigma \varepsilon$, where
$m(\bX)$ is generated from $\exp \{-\frac
{1}{2}(X_1^2/0.8^2+X_2^2/0.9^2+X^2_{3}/1.0+X_4^2/1.1^2) \}$ by
appropriately scaling it to have zero mean and unit variance,
predictors are independently generated from standard normal
distribution and $\varepsilon\sim N(0,1)$ is independent of $X_j$'s.
We set $p=1000$ and $n=100$, and vary the noise level as 0.5 and 1.0,
respectively. We apply the proposed screening method to reduce the
number of predictors from 1000 to 20 and compare it with the method of
\citet{FangFengSongJASA2011} and additional two methods in \citet{Zhuetal2011JASA} and \citet{Lietal2012JASA}.
In this example, the signals are strongest locally at 0 while decay
exponentially fast at other locations, and $X_1$ is the strongest and
$X_4$ is the weakest in their signal strength according to the
coefficients. Note that the iterative screening of \citet{FangFengSongJASA2011} is residual-based while that of \citet{Zhuetal2011JASA} is
projection-based. Thus to be fair, we only compare in terms of the
noniterative version. The frequencies of important predictors being
\begin{table}
\caption{Simulation results for Example~\protect\ref{exa4}}\label{tbex3}
\begin{tabular*}{\textwidth}{@{\extracolsep{\fill}}lccccc@{}}
\hline
$\bolds{\sigma}$ & \textbf{Method} & $\bolds{X_1}$ & $\bolds{X_2}$ & $\bolds{X_3}$ & $\bolds{X_4}$\\
\hline
0.5 &EL & 96 & 92 & 81 & 63 \\
&\citet{FangFengSongJASA2011} & 92& 77 & 62 & 32 \\
& \citet{Zhuetal2011JASA}&\phantom{0}1 & \phantom{0}1 & \phantom{0}5 & \phantom{0}3\\
&\citet{Lietal2012JASA} & 58& 29 & 26 & \phantom{0}9\\
&\citeauthor{ChangTangWu2013} (\citeyear{ChangTangWu2013}) &\phantom{0}0 & \phantom{0}1 & \phantom{0}2 & \phantom{0}5\\[3pt]
1.0 &EL & 81 & 69 & 60 & 36 \\
&\citet{FangFengSongJASA2011} & 73 & 62& 45 & 19 \\
& \citet{Zhuetal2011JASA}&\phantom{0}1 & \phantom{0}0 & \phantom{0}7 & \phantom{0}4\\
&\citet{Lietal2012JASA} & 37 & 12 & 17 & \phantom{0}7\\
&\citeauthor{ChangTangWu2013} (\citeyear{ChangTangWu2013}) &\phantom{0}1 & \phantom{0}1 & \phantom{0}2 & \phantom{0}4\\
\hline
\end{tabular*}\vspace*{3pt}
\end{table}
selected over 100 repetitions are reported in Table~\ref{tbex3} for
different methods. We see that our method performs much better than
that of \citet{FangFengSongJASA2011} thanks to the merit of our method in
detecting local contributions. Additionally, we see that correlation
based methods completely fail in this case while the distance
correlation based method of \citet{Lietal2012JASA} can still detect
signal, while our method performs the best.
\end{exa}

\begin{exa}\label{exa5}
We consider the varying coefficient model in this
example. We generate data from model $ Y=X_1 \beta_1(Z)+X_2\beta_2(Z) +
X_3\beta(Z) + X_4\beta(Z) + \varepsilon$, where predictors are
multivariate normal with $\mathbb{E}(X_j)=0$,\break $\operatorname{var}(X_j)=1$, and
zero correlation,
$\varepsilon\sim N(0,0.1)$ is independent of $X_j$'s, and $Z$ is
independently generated from the standard uniform distribution over
$[0,1]$. The varying
coefficients are given by $\beta_1(z)=\sin(2\pi z+\frac{\pi}{4})$,
$\beta_2(z)=\sin(2\pi z)$, $\beta_3(z)=\cos(2\pi z)$ and
$\beta_4(z)=\sin(2\pi z+\frac{3\pi}{4})$.
We fix the dimensionality $p=1000$ and vary the sample size from 100 to
200. We try to reduce the dimensionality from 1000 to 20 and compare
\begin{table}[t]
\caption{Simulation results for Example~\protect\ref{exa5}}\label{tbex4}
\begin{tabular*}{\textwidth}{@{\extracolsep{\fill}}lccccc@{}}
 \hline
$\bolds{n}$& \textbf{Method} & $\bolds{X_1}$ & $\bolds{X_2}$ & $\bolds{X_3}$ & $\bolds{X_4}$\\
\hline
100     & EL                    &\phantom{0}97& \phantom{0}93& \phantom{0}96& \phantom{0}96\\
                        &\citet{Songetal2014}& \phantom{0}84& \phantom{0}85& \phantom{0}82& \phantom{0}89\\
                        &\citet{Liuetal2014}  &\phantom{0}92& \phantom{0}98& \phantom{0}88& \phantom{0}98\\
                        &\citet{Fanet2014} &\phantom{0}93& \phantom{0}95& \phantom{0}97& \phantom{0}99\\[3pt]
    &\citet{FangFengSongJASA2011}  & \phantom{0}13 &   \phantom{00}8  &  \phantom{00}9 &  \phantom{0}11\\
    &\citet{Zhuetal2011JASA}& \phantom{00}3  &  \phantom{00}6   & \phantom{00}4  &  \phantom{00}5\\
    &\citet{Lietal2012JASA}  &\phantom{00}4  &  \phantom{00}6   & \phantom{00}8  &  \phantom{00}9\\
    &\citeauthor{ChangTangWu2013} (\citeyear{ChangTangWu2013}) &\phantom{00}4  &  \phantom{00}2   & \phantom{00}4  &  \phantom{00}3 \\[6pt]
200    & EL                    & 100& 100 &100& 100\\
                        &\citet{Songetal2014}& 100& 100 &100& 100\\
                        &\citet{Liuetal2014}  & 100& 100 &100& 100\\
                        &\citet{Fanet2014} & 100& 100 &100& 100\\[3pt]
                          &\citet{FangFengSongJASA2011}  & \phantom{0}16  & \phantom{0}13  & \phantom{0}11 &  \phantom{0}15\\
    &\citet{Zhuetal2011JASA}& \phantom{00}5  &  \phantom{00}9  &  \phantom{00}3 &   \phantom{00}6\\
    &\citet{Lietal2012JASA}  & \phantom{00}8  & \phantom{0}15  &  \phantom{00}7 &   \phantom{00}7 \\
    &\citeauthor{ChangTangWu2013} (\citeyear{ChangTangWu2013})  & \phantom{00}2  &  \phantom{00}7  &  \phantom{00}1 &   \phantom{00}3\\
\hline
 \end{tabular*}\vspace*{3pt}
\end{table}
our method to \citet{Fanet2014}, \citet{Liuetal2014} and
\citet{Songetal2014} in terms of the noniterative version due to the
same reason as in the above example.
Table~\ref{tbex4} summarizes the results over 100 repetitions in terms
of how often important predictors are selected. It shows that our
methods perform competitively. Note that the methods of \citet{Fanet2014}, \citet{Liuetal2014}, \citet{Songetal2014} are developed
specially for the varying coefficient model.
\end{exa}

\section{Discussion}\label{s6}

We have proposed and investigated a local independent feature screening method using the marginal empirical
likelihood in conjunction with marginal kernel smoothing methods to detect contributing explanatory variables in a general model setting.
We show that our method is broadly applicable in a wide class of nonparametric and semiparametric models for high-dimensional data analysis.
Theory and numerical examples show that our approach works promisingly.  When the minimal signal is weak or the collinearity  level among the explanatory variables is high,
independence feature screening methods will face substantial difficulty.  How to solve the variable selection problem under such a scenario remains open, and we hope to work
along this direction with the marginal empirical likelihood approach.

Our method is based on the empirical likelihood, and thus necessarily inherits its intensive computation.
Fortunately, the marginal screening methods are highly scalable
by exploring the response variable's dependence on each individual predictor at a time. Consequently,
they are naturally suited for parallel computing. With parallel computing, the computational intensiveness issue of our new method
can be alleviated significantly, making it a practically appealing candidate method.

\section*{Acknowledgements}
We thank the Co-Editor, the Associate Editor and three referees for very constructive comments and suggestions which
have improved the presentation of the paper.


\begin{supplement}[id=suppA]
\stitle{Supplement to ``Local independence feature screening for nonparametric and
semiparametric models by marginal empirical likelihood''}
\slink[doi]{10.1214/15-AOS1374SUPP}  
\sdatatype{.pdf}
\sfilename{aos1374\_supp.pdf}
\sdescription{This supplement contains a real data analysis and all technical proofs.}
\end{supplement}


\printaddresses
\end{document}